# DISCUSSION: THE DANTZIG SELECTOR: STATISTICAL ESTIMATION WHEN $p$ IS MUCH LARGER THAN $n$

By Bradley Efron[1], Trevor Hastie[2] and Robert Tibshirani[3]

*Stanford University*

**1. Introduction.** This is a fascinating paper on an important topic: the choice of predictor variables in large-scale linear models. A previous paper in these pages attacked the same problem using the "LARS" algorithm (Efron, Hastie, Johnstone and Tibshirani [3]); actually three algorithms including the Lasso as middle case. There are tantalizing similarities between the Dantzig Selector (DS) and the LARS methods, but they are not the same and produce somewhat different models. We explore this relationship with the Lasso and LARS here.

**2. Dantzig selector and the Lasso.** The definition of the Dantzig selector (DS) in (1.7) can be re-expressed as

$$(1) \qquad \min_\beta \|X^T(y - X\beta)\|_{\ell_\infty} \quad \text{subject to} \quad \|\beta\|_{\ell_1} \leq s.$$

This makes it look very similar to the Lasso (Tibshirani [6]), or basis pursuit (Chen, Donoho and Saunders [1]):

$$(2) \qquad \min_\beta \|y - X\beta\|_{\ell_2} \quad \text{subject to} \quad \|\beta\|_{\ell_1} \leq s.$$

With a bound on the $\ell_1$ norm of $\beta$, Lasso minimizes the squared error while DS minimizes the maximum component of the gradient of the squared error function. If $s$ is large so that the constraint has no effect, then these are the same. However, for other values of $s$, they are a little different; see Figure 1.

The least angle regression (LARS) algorithm (Efron, Hastie, Johnstone and Tibshirami [3]) for solving the Lasso path makes them look tantalizingly

Received January 2007.
[1]Suported in part by NSF Grant DMS-00-72360 and by NIH Grant 8R01 EB002784.
[2]Suported in part by NSF Grant DMS-05-05676 and by NIH Grant 2R01 CA 72028-07.
[3]Suported in part by NSF Grant DMS-99-71405 and by NIH Contract N01-HV-28183.

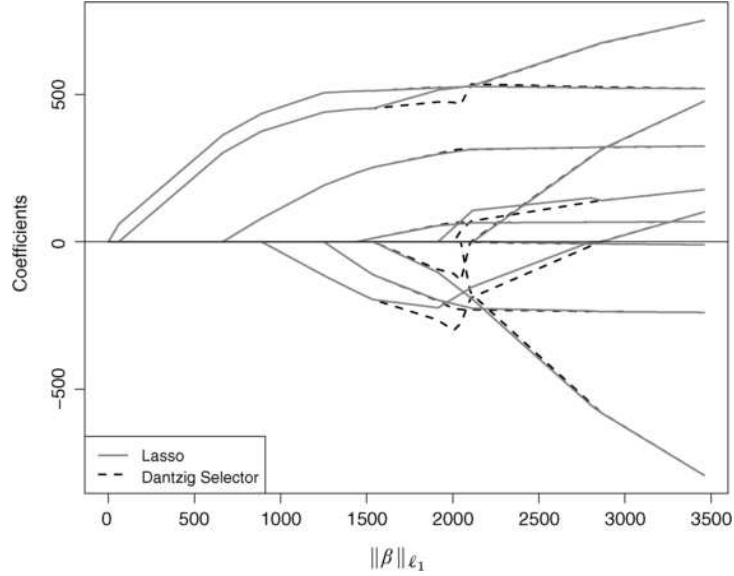

FIG. 1. *The Lasso and DS regularization paths for the diabetes data are mostly identical. The predictors are standardized to have mean zero and unit $\ell_2$ norm, and were used to illustrate the LARS algorithms cited in the text.*

close (see also the homotopy algorithm of Osborne, Presnell and Turlach [4]). In LARS, we start with $\beta = 0$ and identify the predictor having maximal absolute inner product with $y$. We then increase/decrease its coefficient (depending on the sign of the inner product), which in turn reduces its absolute inner product with the current residual $r = y - X\hat{\beta}$. We continue until some other predictor has as large an absolute inner product with the current residual. That predictor is then included in the model, and we move both coefficients in the least squares *equiangular* direction, which keeps their maximal inner products with the residuals the same and decreasing. This process is continued, each time including variables into the model when their inner products catch up with the maximal inner products. Eventually all the inner products are zero, and the algorithm stops. If in addition we drop a predictor out of the model as its coefficient passes through zero, then this LARS algorithm delivers the entire solution set for the Lasso problem (2) for $s$ running from 0 to $\infty$.

Thus at any stage in the Lasso path, the predictors $X_j$ in the model all have equal absolute inner product $|X_j^T(y - X\hat{\beta})|$ with the residuals, and the predictors not in the model have a lower inner product. This is also reflected in the Karush–Kuhn–Tucker conditions for the Lagrange form of (2),

$$\min_{\beta} \tfrac{1}{2}\|y - X\beta\|_{\ell_2}^2 + \lambda\|\beta\|_{\ell_1}, \tag{3}$$

DISCUSSION 3

TABLE 1
*Results for the Lasso and the Dantzig selector on the diabetes data with 64 variables (first 12 shown)*

| Variable $j$ | Lasso | | Dantzig selector | |
|---|---|---|---|---|
| | $X_j^T(y - X\hat{\beta})$ | $\hat{\beta}_j$ | $X_j^T(y - X\hat{\beta})$ | $\hat{\beta}_j$ |
| 1 | 27.4134 | 0.0000 | 26.0046 | 0.0000 |
| 2 | $-83.6413$ | $-77.0062$ | $-83.4945$ | $-73.0993$ |
| 3 | 83.6413 | 502.8695 | 62.5323 | 543.7634 |
| 4 | 83.6413 | 233.5998 | 83.4945 | 223.6250 |
| 5 | $-41.1153$ | 0.0000 | $-43.5949$ | 0.0000 |
| 6 | $-33.8190$ | 0.0000 | $-37.0429$ | 0.0000 |
| 7 | $-83.6413$ | $-164.0632$ | $-83.4945$ | $-155.4648$ |
| 8 | 51.2581 | 0.0000 | 50.5638 | 0.0000 |
| 9 | 83.6413 | 463.4805 | 83.4945 | 455.3289 |
| 10 | 83.6413 | 4.9767 | 83.4945 | 0.0000 |
| 11 | 76.1206 | 0.0000 | 75.6962 | 0.0000 |
| 12 | 83.6413 | 29.7423 | 83.4945 | 13.1410 |
| ⋮ | ⋮ | ⋮ | ⋮ | ⋮ |

The Lasso and DS solutions have the same $\ell_1$ norm $\|\hat{\beta}\|_{\ell_1} = 1734.79$.

which require that

(4) $\qquad X_j^T(y - X\beta) = \lambda \cdot \text{sign}(\beta_j) \qquad \text{for } |\beta_j| > 0,$

(5) $\qquad |X_k^T(y - X\beta)| \leq \lambda \qquad \text{for } |\beta_j| = 0.$

The DS procedure seeks to minimize this maximal inner product directly. How are these different? Table 1 shows an example. The data are the larger version of the diabetes data, consisting of $n = 442$ observations and $p = 64$ predictors (main effects and interactions). The variables have been standardized to have mean zero and variance 1. We have computed both the Lasso and DS solutions with $\|\hat{\beta}\|_{\ell_1} = 1734.79$. At this point, both the Lasso and DS have 12 nonzero coefficients. We give information for the first 12 predictors in the table. We see that *in DS there is a variable (#10) attaining the maximum inner product that is not in the current model*. This is in contrast to the Lasso, where the variables that achieve the maximal inner product are exactly the ones with nonzero coefficients, a consequence of the KKT conditions (4)–(5). DS does this in order to achieve a lower maximal inner product, here 83.49 versus 83.69 for the Lasso. On the other hand, DS gives variable #3 the largest coefficient (actually the largest among all 64 coefficients), while its inner product with the residual is much smaller than that of other variables. As it should, the Lasso solution achieves smaller mean squared error than DS (2827.4 vs. 2829.4).

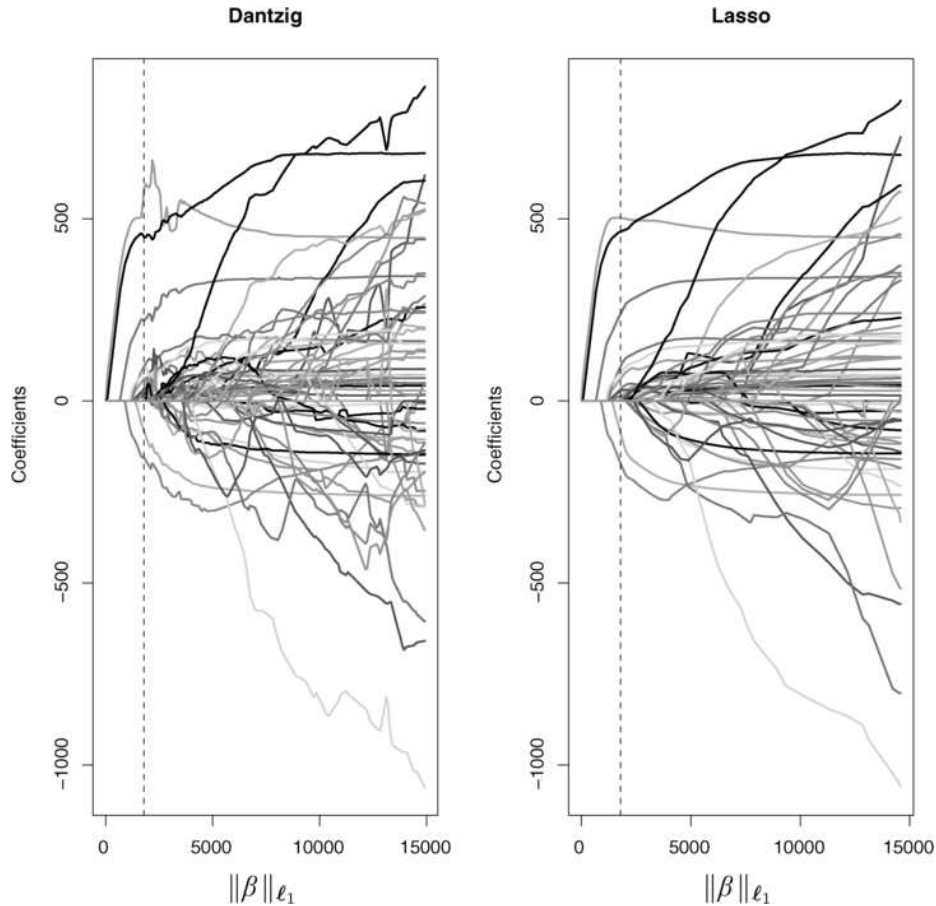

FIG. 2. *Coefficient profiles as a function of $\|\beta\|_{\ell_1}$ for the Dantzig selector and Lasso. There are 64 predictors, the main effects and interactions for the diabetes data. Both paths were truncated at one quarter the norm of the full least squares fit, to allow us to zoom in on the earlier, more relevant parts of the paths.*

We found this surprising and somewhat counterintuitive. In reducing RSS($\beta$) maximally per unit increase in $\|\beta\|_{\ell_1}$, the active set for Lasso does correspond to the variables with largest gradients. We would have also guessed that these gradients were being reduced as fast as possible, but the DS shows this is not the case.

Figure 2 shows the entire solution paths for Lasso and DS for the diabetes data. We see that the DS paths are generally wilder than those of Lasso.

How does this behavior of DS affect its accuracy in practice? We investigate this next.

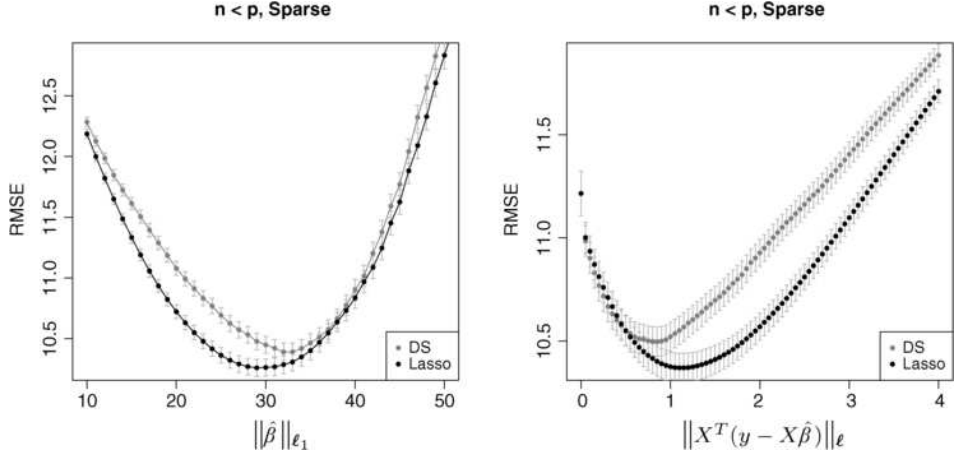

Fig. 3. *RMSE curves for Lasso and DS for the simulation with $n = 15$, $p = 100$, and a sparse coefficient vector $\beta$ with 15 nonzero entries. The left panel uses a grid on $\|\hat{\beta}\|_{\ell_1}$, while the right uses a grid on $\|X^T(y - X\hat{\beta})\|_{\ell_\infty}$.*

**3. Comparison of prediction accuracy.** We conducted a small simulation study to compare the Lasso and DS. We generated data from the model

(6) $$y = X\beta + \varepsilon,$$

with $X$ a matrix of $p = 100$ variables (columns) and $n = 25$ samples. Each of the entries in $X_j$ as well as those in $\varepsilon$ were generated i.i.d. from a Gaussian distribution $N(0,1)$. The first 15 coefficients of $\beta$ were generated from a $N(0,16)$ distribution, and the remainder were set to zero. Hence we dubbed this the $n < p$ *sparse* case. The Lasso and DS coefficient paths were computed on a grid of values for $\|\beta\|_{\ell_1}$, and for each value we computed the root-mean-squared error between $\hat{\beta}$ and the true $\beta$. This was repeated 1000 times, generating a new $X$ and $\varepsilon$ each time, but using the same value for $\beta$. Thus for each value of $\|\hat{\beta}\|_{\ell_1}$, we have 1000 RMSE values corresponding to each of Lasso and DS. Figure 3 (left panel) shows the average and standard deviation for these RMSEs. Lasso is consistently below DS. The right panel shows a similar simulation, except here we use a grid of values for $\lambda = \|X^T(y - X\hat{\beta})\|_{\ell_\infty}$. For DS this amounts to solving the equivalent problem to (1):

(7) $$\min_\beta \|\beta\|_{\ell_1} \quad \text{subject to} \quad \|X^T(y - X\beta)\|_{\ell_\infty} = \lambda,$$

and for Lasso, solving the Lagrange form (3). Whichever way we look at these results, Lasso outperforms DS, and achieves a lower minimum.

We repeated this simulation with everything the same except $\beta$ was *dense*: none was zero and each was generated i.i.d. $N(0,1)$. Figure 4 (left panel) shows the results; here the differences are small. The right panel shows a similar simulation for the $n > p$ case ($n = 100$, $p = 25$) and five nonzero elements in $\beta$. Here the performance of the two procedures is nearly identical.

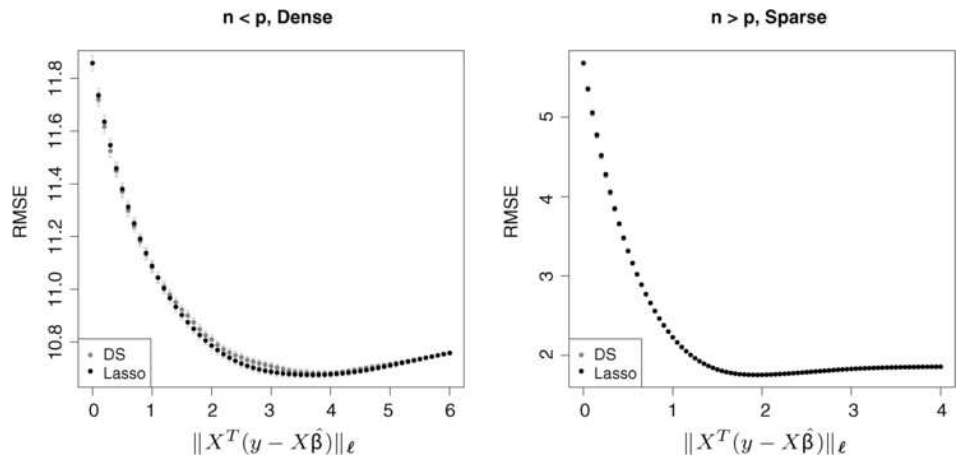

FIG. 4. *RMSE for Lasso and DS in the $n < p$ dense case (left panel) and $n > p$ sparse case (right panel). In the $n < p$ dense case, Lasso slightly outperforms DS but the differences are small. In the $n > p$ sparse case, the performance of the procedures is not distinguishable; this is also the case for the $n > p$ dense case* (*not shown*).

**4. Computational considerations.** The DS problem (1) is a linear program (LP) while the Lasso (2) is a quadratic program. The LARS algorithm for computing the Lasso path is piecewise linear, and the computational load for computing the entire path is equivalent to solving a single least squares problem in the final set of variables. For $n \ll p$ and $\beta$ sparse, Donoho and Tsaig [2] argue that this is the most efficient way to solve any of the Lasso problems. DS will also have a piecewise-linear path algorithm (Rosset and Zhu [5]), but from Figure 2 it is clear that it has many more steps, and is unlikely to provide a similar advantage.

**5. Conclusions.** The optimality properties of the Dantzig selector established by the authors are impressive. We wonder if similar properties hold for the Lasso, and hope that the authors can shed some light on this.

From our brief study, the inherent criterion in DS for including predictors in the model appears to be counterintuitive, and its prediction accuracy seems to be similar to that of the Lasso in some settings, and inferior in other settings. Hence we find little reason to recommend the Dantzig selector over the Lasso.

B. EFRON
T. HASTIE
R. TIBSHIRANI
DEPARTMENT OF STATISTICS AND
 HEALTH RESEARCH AND POLICY
STANFORD UNIVERSITY
STANFORD, CALIFORNIA 94305
USA
E-MAIL: efron@stanford.edu
 hastie@stanford.edu
 tibs@stanford.edu